\documentclass[a4paper,12pt]{article}
\usepackage[T1]{fontenc}
\usepackage[latin1]{inputenc}
\usepackage[english]{babel}
\usepackage{hyperref}
\usepackage[dvips]{graphics}
\usepackage{amsmath}
\usepackage{amssymb}
\usepackage{amsopn}
\usepackage{amsbsy}
\usepackage{amstext}
\usepackage{latexsym}
\usepackage{amsfonts}
\usepackage{array}
\usepackage{epsfig}
\usepackage{color}

\title{How to obtain the continued fraction convergents of the number $e$ by neglecting integrals}
\author{\sc Bakir FARHI \\ {\tt bakir.farhi@gmail.com}}
\date{}

\newtheorem{thm}{Theorem}

\newtheorem{lemma}[thm]{Lemma}

\let\epsilon=\varepsilon
\def\EMts{\mspace{.3mu}}

\def\nb#1{{\left\vert{\EMts\EMts #1 \EMts\EMts}\right\vert}}

\def\EMdash{\leavevmode\hbox to 7.5mm{\vrule height .63ex depth -.59ex
    width 5.4mm\hfill}}

\definecolor{zgroudi}{RGB}{136,77,167}

\begin{document}
\maketitle
\setlength{\parindent}{1cm}
\setlength{\parskip}{0.2cm}

\begin{abstract}
In this note, we show that any continued fraction convergent of
the number $e = 2.71828\dots$ can be derived by approximating some integral
$I_{n , m} := \int_{0}^{1} x^n (1 - x)^m e^x d x$ $(n , m \in
\mathbb{N})$ by $0$. In addition, we present a new way for finding
again the well-known regular continued fraction expansion of $e$.
\end{abstract}
\vspace{1mm}

\section{Introduction, Notations and the Result}~\vspace*{-1cm}

Throughout this paper the expression ``{\it To neglect} a real number'' will mean to approximate that number by $0$.
 
A way of obtaining {\it good} rational approximations of the
number $e$ consists simply to neglect integrals of the form:
\begin{equation}\label{eq0}
I_{n , m} := \int_{0}^{1} x^n (1 - x)^m e^x d x ~~~~~~ (n , m \in
\mathbb{N} ,~ \text{sufficiently large})
\end{equation}
Actually, the neglect of such integrals is justified by the fact
that $I_{n , m}$ tends to $0$ as $n$ and $m$ tend to infinity.
Indeed, by leaning on the Euler $\beta$-function, we have for all
$n , m \in \mathbb{N}$:
$$I_{n , m} \leq e \int_{0}^{1} x^n (1 - x)^m d x = e \cdot \beta(n + 1 , m + 1) = \frac{e}{(n + m + 1) \binom{n + m}{n}} ,$$
which tends to $0$ as $n , m$ tend to infinity. Since $I_{n , m}
\geq 0$, the claimed fact that $I_{n , m}$ tends to $0$ as $n , m$
tend to infinity follows.

\noindent For example, the calculation of $I_{2 , 2}$ gives: $I_{2
, 2} = 14 e - 38$. So if we neglect $I_{2 , 2}$, we obtain the
approximation $e \simeq \frac{19}{7}$ which is a {\it good}
rational approximation of $e$ since it is one of the {\it
convergents} of its regular continued fraction expansion.

The purpose of this paper is essentially to show that any
convergent of the regular continued fraction expansion of the
number $e$ can be obtained by neglecting some integral $I_{n ,
m}$. In addition, we present the paper in a way that the
well-known regular continued fraction expansion of $e$, which is
discovered for the first time by Euler (see e.g., \cite{hw}) and
given by:
\begin{equation}\label{eq1}
e = [2 , 1 , 2 , 1 , 1 , 4 , 1 , \dots] = [2 , \{1 , 2 n , 1\}_{n
\geq 1}]
\end{equation}
will be proved again. Let us define:
\begin{equation}\label{eq2}
e' := [2 , \{1 , 2 n , 1\}_{n \geq 1}] .
\end{equation}
At the end of the paper, we show that $e' = e$ which provides a
new proof of (\ref{eq1}). Our main result is the
following:~\vspace{2mm}

\noindent{\bf Theorem (the main theorem)} {\it Let $n$ be a
natural number and $m$ be a positive integer such that $\nb{n - m}
\leq 1$. Then the neglect of the integral $I_{n , m} :=
\int_{0}^{1} x^n (1 - x)^m e^x d x$ is equivalent to approximate
the number $e$ by some convergent of its continued fraction
expansion.\\ Reciprocally, any convergent of the regular continued
fraction expansion of $e$ can be obtained by neglecting some
integral of the form $I_{n , m}$ with $n , m \in \mathbb{N}$,
satisfying $\nb{n - m} \leq 1$.}~\vspace{2mm}

Now, we are going to give the result which details this theorem.\\
For all positive integer $n$, let $\frac{p_n}{q_n}$ (with $p_n ,
q_n \in \mathbb{N}$, $\gcd(p_n , q_n) = 1$) denotes the
$n$\textsuperscript{th} convergent of the continued fraction
expansion of $e'$. From (\ref{eq2}), we have $e' := [a_1 , a_2 ,
\dots]$, where $a_1 = 2$ and
\begin{equation}\label{eq3}
\begin{cases}
a_{3 k} = 2 k \\
a_{3 k + 1} = a_{3 k - 1} = 1
\end{cases}~~~~~~ (\forall k \geq 1),
\end{equation}
Next, according to the elementary properties of regular continued
fraction expansions (see e.g., \cite{hw}), we have:
\begin{equation}\label{eq4}
\begin{cases}
p_n = a_n p_{n - 1} + p_{n - 2} \\
q_n = a_n q_{n - 1} + q_{n - 2}
\end{cases}~~~~~~ (\forall n \geq 3).
\end{equation}
The details of the main theorem are given by the following:
\begin{thm}[detailing the main theorem]\label{t2}
We have:
\begin{eqnarray}
I_{k , k} & = & (-1)^k k! \left(q_{3 k - 1} e - p_{3 k -1}\right)
~~~~~~ (\forall k \geq 1) \label{eq5} \\
I_{k , k + 1} & = & (-1)^k k! \left(q_{3 k + 1} e - p_{3 k + 1}\right) ~~~~~~ (\forall k \in \mathbb{N}) \label{eq6} \\
I_{k + 1 , k} & = & (-1)^{k + 1} k! \left(q_{3 k} e - p_{3
k}\right) ~~~~~~~~~ (\forall k \geq 1). \label{eq7}
\end{eqnarray}
\end{thm}

\section{The Proofs}~\vspace*{-1cm}

To prove Theorem \ref{t2}, we need the following lemma:
\begin{lemma}\label{lemma}
For all positive integers $n , m$, we have:
\begin{eqnarray}
I_{n , m} & = & m I_{n , m - 1} - n I_{n - 1 , m} \label{eq8} \\
I_{n - 1 , m - 1} & = & I_{n , m - 1} + I_{n - 1 , m} .
\label{eq9}
\end{eqnarray}
\end{lemma}

\noindent{\bf Proof.} Let us prove (\ref{eq8}). By integring by
parts, we have for all positive integers $n , m$:
\begin{eqnarray*}
I_{n , m} & := & \int_{0}^{1} x^n (1 - x)^m e^x d x \\
& = & \int_{0}^{1} x^n (1 - x)^m (e^x)' d x \\
& = & \left[x^n (1 - x)^m e^x\right]_{0}^{1} -
\int_{0}^{1}\left(x^n (1 - x)^m\right)' e^x d x \\
& = & 0 - \int_{0}^{1}\left\{n x^{n - 1} (1 - x)^m - m (1 - x)^{m
- 1} x^n\right\} e^x d x \\
& = & m \int_{0}^{1} x^n (1 - x)^{m - 1} e^x d x - n \int_{0}^{1}
x^{n - 1} (1 - x)^m e^x d x \\
& = & m I_{n , m - 1} - n I_{n - 1 , m} ,
\end{eqnarray*}
as required. Now, let us prove (\ref{eq9}). For all positive
integers $n , m$, we have:
\begin{eqnarray*}
I_{n , m - 1} + I_{n - 1 , m} & = & \int_{0}^{1} x^n (1 - x)^{m -
1} e^x d x + \int_{0}^{1} x^{n - 1} (1 - x)^m e^x d x \\
& = & \int_{0}^{1} x^{n - 1} (1 - x)^{m - 1} \left\{x + (1 -
x)\right\} e^x d x \\
& = & \int_{0}^{1} x^{n - 1} (1 - x)^{m - 1} e^x d x \\
& = & I_{n - 1 , m - 1} ,
\end{eqnarray*}
as required. The lemma is proved.\hfill$\blacksquare$\vspace{3mm}

Now, we are ready to prove Theorem \ref{t2}.\vspace{2mm}

\noindent{\bf Proof of Theorem \ref{t2}.} We proceed by induction
on $k \in \mathbb{N}$. For $k = 0 , 1$, we easily check the
validity of the relations of Theorem \ref{t2}. Now, suppose that
the relations (\ref{eq5}), (\ref{eq6}) and (\ref{eq7}) of Theorem
\ref{t2} hold for some $k \geq 1$ and let us show that they also
hold for the integer $(k + 1)$. Using Formula (\ref{eq8}) of Lemma
\ref{lemma} together with the induction hypothesis, we have:
\begin{eqnarray*}
I_{k + 1 , k + 1} & = & (k + 1) \left(I_{k + 1 , k} - I_{k , k +
1}\right) \\
& = & (k + 1) \left\{(-1)^{k + 1} k! (q_{3 k} e - p_{3 k}) -
(-1)^k k! (q_{3 k + 1} e - p_{3 k + 1})\right\} \\
& = & (-1)^{k + 1} (k + 1)! \left\{(q_{3 k} + q_{3 k + 1}) e -
(p_{3 k} + p_{3 k + 1})\right\} .
\end{eqnarray*}
But according to (\ref{eq3}) and (\ref{eq4}), we have:
$$p_{3 k} + p_{3 k + 1} = p_{3 k + 2} ~~\text{and}~~ q_{3 k} + q_{3 k + 1} = q_{3 k + 2} .$$
Hence:
$$I_{k + 1 , k + 1} = (-1)^{k + 1} (k + 1)! \left(q_{3 k + 2} e - p_{3 k + 2}\right) ,$$
which confirms the validity of Relation (\ref{eq5}) of Theorem
\ref{t2} for $(k + 1)$.

\noindent Next, by using Formulas (\ref{eq8}) and (\ref{eq9}) of
Lemma \ref{lemma} together with the induction hypothesis, we have:
\begin{eqnarray}
I_{k + 1 , k + 2} & = & (k + 2) I_{k + 1 , k + 1} - (k + 1) I_{k ,
k + 2} \notag \\
& = & (k + 2) I_{k + 1 , k + 1} - (k + 1) \left(I_{k , k + 1} -
I_{k + 1 , k + 1}\right) \notag \\
& = & (2 k + 3) I_{k + 1 , k + 1} - (k + 1) I_{k , k + 1} \notag \\
& = & (2 k + 3) (k + 1) (I_{k + 1 , k} - I_{k , k + 1}) - (k + 1)
I_{k , k + 1} \notag \\
& = & (k + 1) \left\{\!\!\!\!\!\!\!\!\!\phantom{\sum}(2 k + 3)
I_{k + 1 , k} - (2 k + 4) I_{k , k
+ 1}\right\} \notag \\
& = & (k + 1) \left\{\!\!\!\!\!\!\!\!\!\phantom{\sum}(2 k + 3)
(-1)^{k + 1} k! (q_{3 k} e - p_{3
k})\right. \notag \\
& ~ & ~~~~~~~~~~~~\left.- (2 k + 4) (-1)^k k! (q_{3 k + 1} e - p_{3 k + 1})\phantom{\sum}\!\!\!\!\!\!\!\!\!\right\} \notag \\
& = & (-1)^{k + 1} (k + 1)!\left\{\!\!\!\!\!\!\!\!\!\phantom{\sum}(2 k + 3) (q_{3 k} e - p_{3 k})
+ (2 k + 4) (q_{3 k + 1} e - p_{3 k + 1})\right\} \notag \\
& = & (-1)^{k + 1} (k +
1)!\left\{\!\!\!\!\!\!\!\!\!\phantom{\sum}\left((2 k + 3) q_{3 k}
+ (2 k + 4) q_{3 k + 1}\right) e \right. \notag \\
& ~ & ~~~~~~~~~~~~~~~~~~~~~~~~\left.- \left((2 k + 3) p_{3 k} + (2
k + 4) p_{3 k + 1}\right)\phantom{\sum}\!\!\!\!\!\!\!\!\!\right\}
\label{eq10}
\end{eqnarray}
But, according to (\ref{eq3}) and (\ref{eq4}), we have:
\begin{eqnarray*}
p_{3 k + 4} & = & p_{3 k + 3} + p_{3 k + 2} \\
& = & ((2 k + 2) p_{3 k + 2} + p_{3 k + 1}) + (p_{3 k
+ 1} + p_{3 k}) \\
& = & (2 k + 2) p_{3 k + 2} + 2 p_{3 k + 1} + p_{3 k} \\
& = & (2 k + 2) (p_{3 k + 1} + p_{3 k}) + 2 p_{3 k + 1} + p_{3 k}
\\
& = & (2 k + 3) p_{3 k} + (2 k + 4) p_{3 k + 1}
\end{eqnarray*}
and similarly, we get:
$$q_{3 k + 4} = (2 k + 3) q_{3 k} + (2 k + 4) q_{3 k + 1} .$$
It follows from (\ref{eq10}) that:
$$I_{k + 1 , k + 2} = (-1)^{k + 1} (k + 1)! (q_{3 k + 4} e - p_{3 k + 4}) ,$$
which confirms the validity of Relation (\ref{eq6}) of Teorem
\ref{t2} for the integer $(k + 1)$.

\noindent Finally, by still using the formulas of Lemma
\ref{lemma} together with the formulas for $I_{k + 1 , k + 1}$ and
$I_{k + 1 , k + 2}$ which we just proved above, we have:
\begin{eqnarray*}
I_{k + 2 , k + 1} & = & I_{k + 1 , k + 1} - I_{k + 1 , k + 2} \\
& = & (-1)^{k + 1} (k + 1)! (q_{3 k + 2} e - p_{3 k + 2}) -
(-1)^{k + 1} (k + 1)! (q_{3 k + 4} e - p_{3 k + 4}) \\
& = & (-1)^{k + 1} (k + 1)!
\left\{\!\!\!\!\!\!\!\!\!\phantom{\sum}(q_{3 k + 2} - q_{3 k + 4})
e - (p_{3 k + 2} - p_{3 k + 4})\right\} .
\end{eqnarray*}
But since (according to (\ref{eq3}) and (\ref{eq4})): $q_{3 k + 4}
= q_{3 k + 3} + q_{3 k + 2}$ and $p_{3 k + 4} = p_{3 k + 3} + p_{3
k + 2}$ then it follows that:
$$I_{k + 2 , k + 1} = (-1)^{k + 2} (k + 1)! (q_{3 k + 3} e - p_{3 k + 3}) ,$$
which confirms the validity of Relation (\ref{eq7}) of Theorem
\ref{t2} for the integer $(k + 1)$.\\
The three relations of Theorem \ref{t2} thus hold for $(k + 1)$.
The proof of Theorem \ref{t2} is
complete.\hfill$\blacksquare$\vspace{2mm}

Theorem \ref{t2} permits us to establish a new proof for the fact
that the regular continued fraction expansion for the number $e$
is given by (\ref{eq1}).~\vspace{2mm}

\noindent{\bf A new Proof of (\ref{eq1}).} Relation (\ref{eq5}) of
Theorem \ref{t2} shows that for all positive integer $k$, we have:
\begin{eqnarray}
\nb{e - \frac{p_{3 k - 1}}{q_{3 k - 1}}} & = & \frac{\nb{I_{k ,
k}}}{k! q_{3 k - 1}} \notag \\
& \leq & \frac{I_{k , k}}{k!} ~~~~~~ (\text{since $I_{k , k} \geq
0$ and $q_{3 k - 1} \in \mathbb{Z}_+^*$}) \label{eq11}
\end{eqnarray}
Next, by using the simple inequalities $x (1 - x) \leq
\frac{1}{4}$ and $e^x \leq e$ ($\forall x \in [0 , 1]$), we have
for all positive integer $k$:
$$I_{k , k} := \int_{0}^{1}\left(x (1 - x)\right)^k e^x d x \leq \int_{0}^{1} \frac{e}{4^k} d x = \frac{e}{4^k} .$$
It follows by inserting this in (\ref{eq11}) that:
$$\nb{e - \frac{p_{3 k - 1}}{q_{3 k - 1}}} \leq \frac{e}{4^k k!} ~~~~~~ (\forall k \geq 1) ,$$
which shows that $\frac{p_{3 k - 1}}{q_{3 k - 1}}$ tends to $e$ as
$k$ tends to infinity. But since $\frac{p_{3 k - 1}}{q_{3 k - 1}}$
represents the $(3 k - 1)$\textsuperscript{th} convergent of the
regular continued fraction expansion of $e'$, we have on the other
hand $\lim_{k \rightarrow + \infty} \frac{p_{3 k - 1}}{q_{3 k -
1}} = e'$. Hence $e = e'$, which confirms
(\ref{eq1}).\hfill$\blacksquare$~\vspace{2mm}

\section{Remarks about the analog of the main \\ theorem concerning
the number $\pi$:}~\vspace*{-1cm}

The analogs of the integrals $I_{n , m}$ whose the neglect leads
to approximate the number $\pi$ by the convergents of its regular
continued fraction expansion are not known in their general form.
However, for the particular famous approximation $\pi \simeq
\frac{22}{7}$, Dalzell \cite{dal} noticed that
$$\int_{0}^{1} \frac{x^4 (1 - x)^4}{1 + x^2} d x = \frac{22}{7} - \pi .$$
So the neglect of the later integral leads to the
Archimedes approximation $\pi \simeq \frac{22}{7}$.\\
For some other continued fraction convergents of $\pi$ (like
$\frac{333}{106} , \frac{355}{113} , \dots$), Lucas \cite{lucas}
experimentally obtained some integrals having the form:
$$\int_{0}^{1} \frac{x^n (1 - x)^m (a + b x + c x^2)}{1 + x^2} d x ~~~~~~ (n , m , a , b , c \in \mathbb{N})$$
whose the neglect leads to approximate $\pi$ by those convergents.

~\\


\end{document}